\documentclass[twoside,leqno]{article}
\usepackage{latexsym}
\usepackage{amssymb}
\input{epsf.tex}
\newcommand{\uno}{\hbox{\bf 1}}
\newcommand{\Var}{\mathop{\rm Var}}
\newcommand{\barF}{\overline{F}}
\newcommand{\Prob}{\mathop{\rm P}}
\newcommand{\remove}[1]{}
\newcommand{\Reali}{\mathbb{R}}
\title{\normalsize \bf ON CERTAIN BOUNDS FOR FIRST-CROSSING-TIME PROBABILITIES \\ 
OF A JUMP-DIFFUSION PROCESS\footnote{Paper appeared in Sci.\ Math.\ Jpn.\ 64 (2006), no. 2, 449--460.}}
\author{\normalsize \sc Antonio Di Crescenzo,\quad Elvira Di Nardo \ and \ Luigi M.\ Ricciardi}
\date{\empty}
\begin{document}
\maketitle
\begin{abstract}
We consider the first-crossing-time problem through a constant boundary
for a Wiener process perturbed by random jumps driven by a counting process. 
On the base of a sample-path analysis of the jump-diffusion process
we obtain explicit lower bounds for the first-crossing-time density 
and for the first-crossing-time distribution function. In the case of the 
distribution function, the bound is improved by use of processes comparison 
based on the usual stochastic order. The special case of constant jumps 
driven by a Poisson process is thoroughly discussed.

\medskip\noindent
{\em AMS Classification:} 
60G40,  
60J65,  
60E15   

\medskip\noindent
{\em Key words and phrases:} Jump-diffusion process; Wiener process; first-crossing time;
usual stochastic order.

\end{abstract}
%
\section{Introduction} 
In a variety of applied contexts a relevant role is played by jump-diffusion processes, 
i.e.\ by diffusion processes to which jumps occurring at random times are superimposed. 
Indeed, such processes are for example invoked for the description of stochastic 
neuronal activity (see Giraudo and Sacerdote \cite{GirSac97} and Giraudo {\em et al.}\ 
\cite{GiSaSi02}), of complex queueing systems (see Perry and Stadje \cite{PeSt01}), 
of random assets in mathematical finance (see Ball and Roma \cite{BaRo93}), 
of surplus of insurance companies in ruin theory (see Gerber and Landry \cite{GeLa98}), 
of acto-myosin interaction in biomathematics (see Buonocore {\em et al.}\ 
\cite{Buonocore}).  Despite the relevance of the first-crossing-time (FCT) 
problem for jump-diffusion processes in such contexts, only few analytical 
results are available on their probability density function (pdf), even in 
the case of very simple boundaries. The available results are mainly focused on 
equations involving the FCT moments (see, for instance, Abundo~\cite{Abundo}, 
Giraudo and Sacerdote~\cite{GirSac96}, and Tuckwell \cite{Tuckwell}), which however 
seem to be hardly manageable for practical purposes. Other results concerning 
bounds, obtained by use of Laplace transform, are limited to upper and lower 
bounds for the ruin probability in a jump-diffusion process involved in a risk 
model perturbed by Brownian motion (see Yin and Chiu \cite{YinChiu04}), and 
for the mean and the variance of the hitting time in certain jump-diffusion 
processes (see Sch\"{a}l \cite{Sc93}). Analytical results on FCT pdf's are 
also very rare (see Kou and Wang \cite{KoWa03}, where explicit solutions of 
the Laplace transform of the distributions of the first-crossing times are 
disclosed for a Brownian motion perturbed by double exponentially distributed 
jumps). Hence, as a viable alternative, efficient algorithms have been devised 
in order to evaluate FCT densities (cf.\ the recent contributions by 
Atiya and Metwally \cite{AtMe2005}, and by Di Crescenzo {\em et al.}\ \cite{DCDNR2005}).
\par
In order to evaluate the performance of simulation algorithms, one needs to 
come up with some sample cases in which FCT densities and distribution functions 
obtained by simulation can be compared with the corresponding bounds 
analytically determined. Hence, in the present paper lower bounds for FCT 
densities and distribution functions will be determined for jump-diffusion 
processes based on the Wiener process in the presence of a constant boundary. 
Our results are based on a sample-path analysis of the jump-diffusion process 
and on some specific features of the underlying Wiener process, such as the 
space and time homogeneity and the availability of a closed form of the FCT 
density through a constant boundary. It must be pointed out that our approach 
is at all different from that of Bischoff and Hashorva \cite{BiHa2005}, 
where the Cameron-Martin-Girsanov formula is exploited to obtain a lower bound 
for the boundary crossing probability of Brownian bridge with trend. 
\par 
In Section~2 we formally describe the jump-Wiener model, that consists of the 
superposition of a Wiener process and of a jump process with generally distributed 
jumps occurring at the occurrences of a counting process. The FCT problem for 
such process through a constant boundary is then addressed and a lower bound for 
the FCT pdf is then determined in Section~3. Use of such a bound is then 
made for the special case when upward and downward jumps have constant amplitudes 
and occur according to a Poisson process. Section~4 presents a lower bound for 
the FCT cumulative distribution function (cdf). The special case when only upward 
constant jumps are allowed, and are separed by random times having exponential 
distribution is thoroughly investigated. In this case the bound is improved by 
making use of a technique based on the comparison of FCT's by the 
``usual stochastic order''. Finally, in Section~5 some remarks on the computational 
aspects are given.  
\section{FCT problem for a jump-diffusion process} 
Let $\{X(t)\}_{t \geq 0}$ be a jump-diffusion process defined by
\begin{equation}
X(t)= W(t) + Y(t)
\label{(model)}
\end{equation}
where $\{W(t)\}_{t \geq 0}$ and $\{Y(t)\}_{t \geq 0}$ are
independent stochastic processes and 
\begin{description}
\item[{\it (i)}] $\{W(t)\}_{t \geq 0}$ is a Wiener process
with drift $\mu \in \Reali,$ and variance $\sigma^2 \in (0,+\infty)$
per unit time, starting at $W(0)=x_0;$
\item[{\it (ii)}] $\{Y(t)\}_{t \geq 0}$ is a jump process such that $Y(0)=0$ 
and $Y(t)=\sum_{i=1}^{N(t)} J_i$, $t>0$, with $Y(t)=0$ when $N(t)=0$ and 
$J_1,J_2,\ldots$ real-valued i.i.d.\ r.v.'s such that $\Prob\{J_i=0\}<1$ 
for all $i$. The cdf of $J_i$ will be denoted by $F_J(x)$ and the survival 
function by $\overline{F}_J(x)=1-F_J(x)$. By $\{N(t)\}_{t\geq 0}$ we denote 
a counting process independent of $\{J_1,J_2,\ldots\}$ and characterized by 
i.i.d.\ absolutely continuous positive renewals $R_1,R_2,\ldots$ having 
pdf $f_R(x)$, cdf $F_R(x)$ and survival function $\overline{F}_R(x)=1-F_R(x)$. 
\end{description}
Let 
$$
 f(x,t\,|\,x_0)
 =\frac{\partial}{\partial x}\Prob\left\{X(t)\leq x\,|\,X(0)=x_0\right\},
 \qquad t>0,
$$
be the conditional pdf of $\{X(t)\}$. Since 
$$
 \Prob\left\{ W(t)+\sum_{i=1}^{k} J_i \in {\rm d}x \, \Bigg| \, W(0)=x_0 \right\} 
 =\int_{\Reali} f_W(x - u,t \, | \, x_0) \, {\rm d}F^{[k]}(u) \, {\rm d}x,  
$$
where $F^{[k]}(u)$ denotes the cdf of $\sum_{i=1}^k J_i$ and 
$$
 f_W(x,t\,|\,x_0)=\frac{1}{\sqrt{2 \pi \sigma^2 t}} 
 \exp\left\{-\frac{(x-x_0-\mu t)^2}{2 \sigma^2 t}\right\}
$$
is the pdf of $W(t)$, from the assumptions on $\{X(t)\}$ we have:
\begin{equation}
f(x,t\,|\,x_0) = \sum_{k=0}^{\infty} P[N(t)=k] \int_{\Reali}
f_W(x-u,t \, | \, x_0) \, {\rm d}F^{[k]}(u).
\label{(denstrans)}
\end{equation}
Denoting by $\phi_N(u,t)$ the probability generating function of $N(t)$ and by 
$\psi_J(s)$ the moment generating function (m.g.f.) of $J_i$, the m.g.f.\ of 
$\{X(t)\}$ is given by 
\begin{equation}
 E \left[ e^{s X(t)}\,\Big|\,X(0)=x_0 \right] 
 = \phi_N[\psi_J(s),t]\,\exp \left\{ (x_0+\mu t)s
 + \frac{\sigma^2 t}{2} s^2 \right\},
\label{(genfun0)}
\end{equation}
and  
\begin{eqnarray}
 E[X(t)\,|\,X(0)=x_0] 
 \!\!\! & = & \!\!\! x_0 + \mu t + E(J_1) \, E[N(t)], \label{(media)} \\
 \Var[X(t)\,|\,X(0)=x_0] 
 \!\!\! & = & \!\!\! \sigma^2 t + \Var(J_1)\, E[N(t)] 
 + E^2(J_1)\, \Var[N(t)]. \label{(varianza)}
\end{eqnarray}
\par
Assuming $x_0<S$, hereafter we shall discuss the FCT problem through 
a constant boundary $S$ for the jump-diffusion process defined in (\ref{(model)}). 
We denote by 
\begin{equation}
 T_X=\inf \{t \geq 0: X(t) \geq S\}, 
 \qquad \Prob\{X(0)=x_0\}=1,  
 \label{(defTX)}
\end{equation}
the FCT of $\{X(t)\}$ through $S$ from below, and by 
\begin{equation}
 g_X(S,t \, | \, x_0) = \frac{\partial}{\partial t} 
 \Prob \left\{ T_X \leq t\,|\, X(0) = x_0 \right\}, \qquad t > 0 
 \label{(eq2)}
\end{equation}
its pdf. The corresponding cdf will be denoted by
\begin{equation}
 G_X(S,t\,|\,x_0)=\Prob(T_X\leq t\,|\,X(0)=x_0), \qquad t \geq 0. 
 \label{(eqG)}
\end{equation}
To determine the ultimate FCT probability and the FCT moments one is usually asked 
to solve appropriate integro-differential equations (see Abundo~\cite{Abundo} 
and Tuckwell~\cite{Tuckwell}). The determination of closed-form expressions for 
the FCT pdf and cdf is instead a harder problem, because so far no analytical 
method appear to be available thus for. Hence, in the following sections we 
shall confine our investigation to determining useful lower bounds for 
density (\ref{(eq2)}) and cdf (\ref{(eqG)}). 
\par
Let $T_W$ denote the FCT from below of Wiener process $\{W(t)\}$ from $x_0<S$ 
to the constant boundary $S$. Some well-known results on $T_W$ that will be 
used later are recalled hereafter: 
\newline
$\bullet$ \ FCT cdf of $\{W(t)\}$ through $S$:
\begin{eqnarray}
 G_W(S,t \, | \, x_0)  
 \!\!\!  & = & \!\!\! \Prob\{T_W\leq t\,|\,W(0)=x_0\}  
 \label{(distr)} \\
 & = & \!\!\! \Phi \left(- \frac{S-x_0-\mu t}{\sqrt{\sigma^2 t}} \right) 
 +\exp \left( 2 \mu \frac{S-x_0}{\sigma^2} \right) \,
 \Phi \left(-\frac{S-x_0+\mu t}{\sqrt{\sigma^2 t}} \right),
 \nonumber
\end{eqnarray}
where
\begin{equation}
 \Phi(z)  =  \frac{1}{\sqrt{2 \pi}} \int_{-\infty}^{z}
 e^{- x^2/2} {\rm d}x, 
 \qquad z \in \Reali
 \label{(norm)};
\end{equation}
\newline
$\bullet$ \ FCT pdf of $\{W(t)\}$ through $S$:
\begin{equation}
g_W(S,t \, | \, x_0) = \frac{\partial}{\partial t} \Prob\{T_W\leq t\,|\,W(0)=x_0\}
=  \frac{S-x_0}{\sqrt{2 \pi \sigma^2 t^3}}
\exp \left\{- \frac{(S-x_0-\mu t)^2}{2 \sigma^2 t} \right\}; 
\label{(wiener)}
\end{equation}
\newline
$\bullet$ \ $S$-avoiding transition pdf of $\{W(t)\}$:
\begin{eqnarray}
 \hspace{1cm}\alpha_W(x,t \, | \, x_0) 
 \!\!\! &=&  \!\!\! \frac{\partial}{\partial x} \Prob\{W(t) \leq x, T_W>t\,|\,W(0)=x_0\} 
 \label{(avoid)} \\
 &=& \!\!\! \frac{1}{\sqrt{2 \pi \sigma^2 t}} 
 \,\exp \left\{ - \frac{(x-x_0-\mu t)^2}{2 \sigma^2 t}\right\}
 \left[ 1 \!-\! \exp \left\{ - \frac{2(S-x)(S-x_0)}{\sigma^2 t}\right\}\right]. 
 \nonumber
\end{eqnarray}
%
\section{Lower bound for FCT pdf} 
In order to obtain the preannounced lower bound for pdf (\ref{(eq2)}) let us 
note that, for all $t>0$, event $\{T_X\in (t,t+{\rm d}t)\}$ can be decomposed 
into the following three mutually exclusive events: \\
{\it (i)} \ 
the first jump occurs after $t$ and the first crossing through $S$ occurs in 
$(t,t+{\rm d}t)$ due to the diffusive component of $\{X(t)\}$; \\
{\it (ii)} \ 
in $(0,t)$ the diffusive component of $\{X(t)\}$ does not cross the boundary, 
the first jump occurs in $(t,t+{\rm d}t)$ and it causes the first crossing; \\
{\it (iii)} \ 
the first jump occurs at time $\theta\in(0,t)$ and it does not cause the first 
crossing, the diffusive component of $\{X(t)\}$ having not crossed the boundary 
in $(0,\theta)$; the first crossing finally occurs in $(t,t+{\rm d}t)$. 
\par
Hence, for all $t > 0$ the following equation holds:
\begin{eqnarray}
 g_X(S,t \, | \, x_0) \!\!\! 
 & = & \!\!\! \overline{F}_R(t) \, g_W(S,t \, | \, x_0) 
 + f_R(t) \int_{-\infty}^S \alpha_W(x,t \, | \, x_0) \, \overline{F}_J(S-x)\,{\rm d}x 
 \label{(dens)} \\
 & + & \!\!\! \int_0^t {\rm d}F_R(\theta) \int_{-\infty}^S 
 \alpha_W(x,\theta\,|\,x_0) 
 \int_{-\infty}^{S-x} g_X(S,t-\theta\,|\,x+u)\,{\rm d}F_J(u)\,{\rm d}x.
 \nonumber  
\end{eqnarray}
The formal proof of Eq.\ (\ref{(dens)}) has been given by Giraudo and Sacerdote 
\cite{GirSac02} in the special case of jumps of constant size at the occurrence 
of a Poisson process, aiming to find out conditions under which the 
FCT density becomes multimodal (see also Sacerdote and Sirovich \cite{SaSi03} 
on this topic). As an immediate consequence of (\ref{(dens)}) a lower bound 
for the FCT pdf of $\{X(t)\}$ can be obtained. Indeed, for all $t>0$ the 
following inequality holds:
\begin{eqnarray}
 && g_X(S,t \, | \, x_0) 
 \geq \overline{F}_R(t) \, g_W(S,t \, | \, x_0) 
 + f_R(t) \int_{-\infty}^S \alpha_W(x,t \, | \, x_0) \, \overline{F}_J(S-x)\,{\rm d}x 
 \label{(ricor)} \\
 && \quad + \int_0^t \overline{F}_R(t-\theta) \,{\rm d}F_R(\theta) \int_{-\infty}^S 
 \alpha_W(x,\theta\,|\,x_0) \int_{-\infty}^{S-x} 
 g_W(S,t-\theta\,|\,x+u)\,{\rm d}F_J(u)\,{\rm d}x.
 \nonumber
\end{eqnarray}
Note that a tighter lower bound can be obtained by repeated substitutions of $g_X$ 
in the last term of the right-hand-side of (\ref{(dens)}). However, this would 
include many terms involving progressively high-order integrals, unsuitable for 
computational purposes. 
\subsection{A special case} 
Let us now study a special case of model (\ref{(model)}). First of all, we assume 
that the jumps are separated by i.i.d.\ exponential random times $R_i$, i.e.\ 
restrict our attention to the case when $\{N(t)\}$ is a Poisson process whose 
parameter will be denoted by $\lambda$. In this case pdf $f(x,t\,|\,x_0)$ is 
solution of the following integro-differential equation (see Buonocore {\em et al.}\ 
\cite{Buonocore} or Di~Crescenzo {\em et al.}\ \cite{Martinucci}):
$$ 
 \frac{\partial \, f}{\partial t} 
 = -\lambda f - \mu \frac{\partial \,f}{\partial x} 
 + \frac{\sigma^2}{2} \frac{\partial^2 \, f}{\partial x^2}
 + \lambda \int_{\Reali} f(x-y,t\,|\,x_0) \,{\rm d}F_J(y),
$$
which clearly reduces to the Fokker-Plank equation of a Wiener process 
if $\lambda=0$. Moreover, let us assume that the jumps have positive or 
negative constant amplitude. In other words, the random jumps $J_i$ are 
distributed as 
\begin{equation}
 J_i= \left\{ \begin{array}{ll}
 a & \hbox{w.p. $\eta$} \\
 -b & \hbox{w.p. $1-\eta$,}
 \end{array} \right.
 \label{(vaJ)}
\end{equation}
with $0<\eta<1$, $a>0$ and $b>0$. Hence, $Y(t)$ can be expressed as 
$$
 Y(t) = a\,N_1(t) - b\,N_2(t), 
 \qquad t>0
$$
where $\{N_1(t)\}_{t \geq 0}$ and $\{N_2(t)\}_{t \geq 0}$ are independent 
Poisson processes with rate $\eta\lambda$ and $(1-\eta)\lambda$, respectively. 
From (\ref{(denstrans)}), for $x\in\Reali$ and $t\geq 0$ we obtain
\begin{equation}
 \qquad f(x,t\,|\,x_0) 
 =  \frac{e^{-\lambda t}}{\sqrt{2 \pi \sigma^2 t}}
 \sum_{k=0}^{\infty} \sum_{j=0}^{\infty} 
 \frac{(\eta\lambda t)^j}{j!}\frac{((1-\eta)\lambda t)^k}{k!} 
 \exp \left\{- \frac{(x+ a j - b k -x_0-\mu t)^2}{2 \sigma^2 t} \right\}.
 \label{(densPois)}
\end{equation}
Since under the present assumptions  
$\psi_J(s)= \eta\,e^{a s}+ (1-\eta)\,e^{-b s}$ and 
$\phi_N(u,t)=e^{-\lambda t(1-u)}$, from (\ref{(genfun0)}) the
m.g.f.\ of $\{X(t)\}$ for $s\in\Reali$ becomes
$$
 E\left[e^{s X(t)}\,\Big|\,X(0)=x_0\right] 
 = \exp\left\{-\lambda t\left[1-\eta\,e^{a s}-(1-\eta)\,e^{-b s}\right]
 +(x_0+\mu t)s+ \frac{\sigma^2 t}{2} s^2\right\}.
$$
As already pointed out in Di~Crescenzo {\em et al.}\ \cite{DCDNR2005}, 
from (\ref{(media)}) and (\ref{(varianza)}) we have
\begin{eqnarray*}
 E[X(t)\,|\,X(0)=x_0] \!\!\!
 & = & \!\!\! x_0 + \mu t + [\eta a-(1-\eta) b] \lambda t, \\
 \Var[X(t)\,|\,X(0)=x_0] \!\!\! 
 & = & \!\!\! \sigma^2 t + [\eta a^2 + (1-\eta) b^2] \lambda t.
\end{eqnarray*}
\par
We stress that the lower bound given in (\ref{(ricor)}) can now be explicitly
evaluated. Indeed, under the present assumptions for all $t>0$ we have:
$$
 g_X(S,t\,|\,x_0) \geq g_{\ell}(S,t\,|\,x_0),
$$
with
\begin{eqnarray}
 g_{\ell}(S,t\,|\,x_0) \!\!\!
 & := & \!\!\! e^{-\lambda t} \Bigg\{ g_W(S,t\,|\,x_0) 
 + \eta\lambda \int_{S-a}^S \alpha_W(x,t\,|\,x_0)\,{\rm d}x   
 \label{(bound)} \\
 & + & \!\!\! \eta \lambda\int_0^t {\rm d}\theta \int_{-\infty}^{S -a}
 \alpha_W(x,\theta\,|\,x_0) \, g_W(S,t-\theta\,|\,x+a)\,{\rm d}x
 \nonumber \\
 & + & \!\!\! (1-\eta)\lambda \int_0^t {\rm d}\theta \int_{-\infty}^S 
 \alpha_W(x,\theta\,|\,x_0) \, g_W(S,t-\theta\,|\,x-b)\,{\rm d}x \Bigg\},
 \nonumber 
\end{eqnarray}
where $g_W(S,t\,|\,x_0)$ is given in (\ref{(wiener)}) and 
\begin{eqnarray*}
 && \hspace{-1cm}
 \int_{S-a}^S \alpha_W(x,t\,|\,x_0)\,{\rm d}x 
 =\Phi \left(\frac{S-x_0-\mu t}{\sigma\sqrt{t}} \right) 
 - \Phi \left(\frac{S-a-x_0-\mu t}{\sigma\sqrt{t}} \right)   \\
 && \qquad\qquad -\exp\left\{\frac{2 \mu}{\sigma^2} (S-x_0) \right\}
 \left[\Phi\left(-\frac{S-x_0+\mu t}{\sigma\sqrt{t}}\right)
 - \Phi\left(-\frac{S+a-x_0+\mu t}{\sigma\sqrt{t}}\right)\right],
\end{eqnarray*}
with $\Phi(z)$ defined in (\ref{(norm)}). The right-hand-side of 
(\ref{(bound)}) can be numerically evaluated. Indeed, making use of identity 
$$
 \int_0^{\infty} (z-\delta) e^{-(a z^2 + b z + \gamma)} {\rm d}z  
 = \frac{e^{-\gamma}}{2 a} \left\{ 1 - 2 \sqrt{a \pi} 
 \left(\frac{b}{2a} + \delta \right) e^{\frac{b^2}{4 a}}
 \left[1-\Phi\left(\frac{b}{\sqrt{2 a}} \right)\right] \right\},
 \qquad a >0,
$$
following from Eq.\ 7.4.2 of Abramowitz and Stegun \cite{Abr}, one obtains:
\begin{eqnarray*}
 \int_{-\infty}^{S-m} \alpha_W(x,\theta\,|\,x_0) \, 
 g_W(S,t-\theta\,|\,x+u)\, {\rm d}x \!\!\!
 & = & \!\!\! \frac{1}{2 \pi t} \sqrt{\frac{\theta}{t-\theta}}
 \exp\left\{ - \frac{(S-u-x_0-\mu t)^2}{2 \sigma^2 t} \right\} \\
 & \times & \!\!\! \left\{A_{-}(\theta,u)
 -\exp\left[-\frac{2(S-x_0)u}{\sigma^2 t}\right] A_+(\theta,u) \right\},
\end{eqnarray*}
where $m=\max\{u,0\}$ and 
\begin{eqnarray*}
 A_{\pm}(\theta,u) \!\!\! 
 & = & \!\!\! \exp \left\{- \frac{[(m-u)\theta \pm(S-x_0 \pm m)(t-\theta)]^2}
 {2 \sigma^2 \theta t (t-\theta)} \right\}  \\
 & \mp & \!\!\! \sqrt{\frac{2 \pi (t-\theta)}{\sigma^2 t \theta}}\, 
 (S-x_0 \pm u) \left\{1-\Phi\left[ \frac{(m-u)\theta \pm (S-x_0 \pm m)(t-\theta)}
 {\sqrt{\sigma^2 \theta t (t-\theta)}} \right]\right\}.
\end{eqnarray*}
%
%
\begin{figure}[t]
\begin{center}
\leavevmode
\epsfxsize=300pt
\epsffile{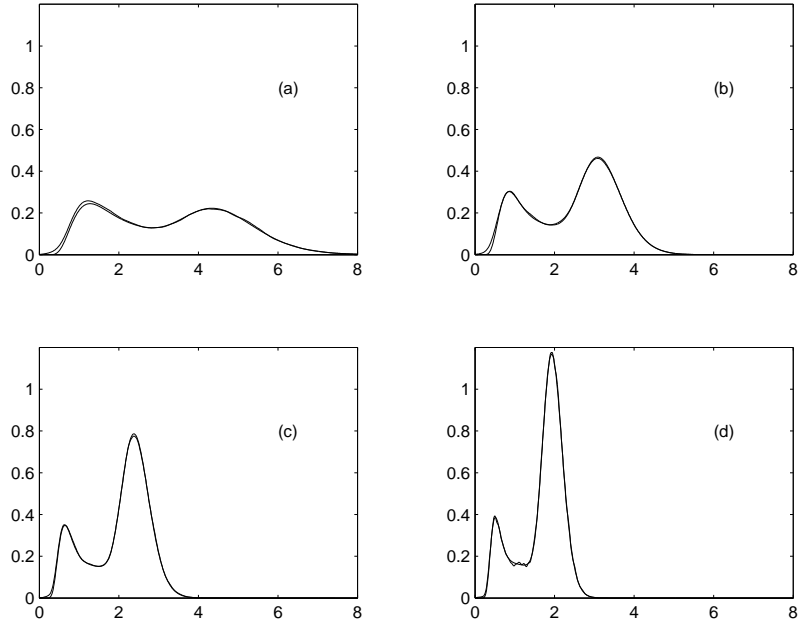}
\caption{
Estimated FCT density $\widehat{g}_X(S,t\,|\,x_0)$ for the Wiener process 
with constant Poisson-paced jumps for $a=b=3.75$, and lower bound 
$g_{\ell}(S,t\,|\,x_0)$, with $x_0=0$, $S=5$, $\lambda=0.2$, 
$\eta=1$, $\sigma^2=0.2$ and (a) $\mu=1.0$, (b) $\mu=1.5$, (c) $\mu=2.0$, (d) $\mu=2.5$.
}
\end{center}
\end{figure}
%
\begin{figure}[t]
\begin{center}
\leavevmode
\epsfxsize=300pt
\epsffile{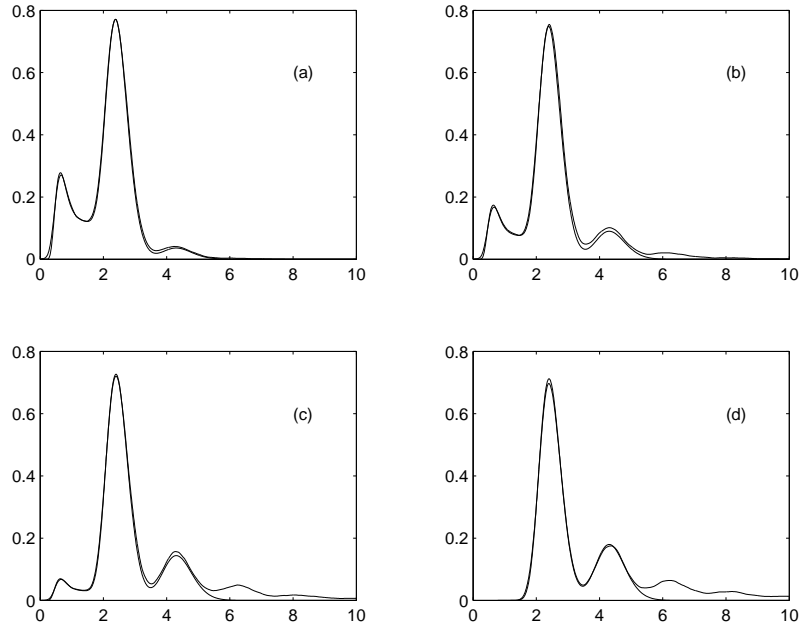}
\caption{
Estimated FCT density $\widehat{g}_X(S,t\,|\,x_0)$ for the Wiener process 
with constant Poisson-paced jumps, and lower bound $g_{\ell}(S,t\,|\,x_0)$, 
for $\mu=2.0$ and (a) $\eta=0.8$, (b) $\eta=0.5$, (c) $\eta=0.2$, (d) $\eta=0$. 
All other parameters are chosen as in Figure 1. 
}
\end{center}
\end{figure}
\par
Figures 1 and 2 show certain multimodal estimates of $g_X(S,t\,|\,x_0)$ obtained by means 
of the simulation procedure described in Di~Crescenzo {\em et al.}\ \cite{DCDNR2005}, 
together with the lower bound given in (\ref{(bound)}). The latter appears to be 
very close to the estimated pdf in various cases; the goodness of the approximation 
is discussed in Section 5. 
\section{Lower bound for FCT cdf} 
In this section we obtain a lower bound for the FCT distribution function (\ref{(eqG)}). 
Similarly to the case of the FCT density we note that event $\{T_X\leq t\}$, $t>0$, 
can be decomposed into three mutually exclusive events: \\
{\it (i)} \ 
the first jump occurs after $t$ and the first crossing through $S$ occurs in 
$(0,t]$ due to the diffusive component of $\{X(t)\}$; \\
{\it (ii)} \ 
the first jump occurs at time $\theta\in(0,t)$ and it causes the first crossing 
through $S$; \\
{\it (iii)} \ 
the first jump occurs at time $\theta\in(0,t)$ and it does not cause the first 
crossing, the diffusive component of $\{X(t)\}$ having not crossed the boundary 
in $(0,\theta)$; the first crossing finally occurs in $(\theta,t]$. 
\par
The following equation thus holds for all $t>0$: 
\begin{eqnarray}
 && \hspace{-0.5cm}
    G_X(S,t\,|\,x_0)=\barF_R(t)\,G_W(S,t\,|\,x_0)
    +\int_0^t {\rm d}F_R(\vartheta)\int_{-\infty}^S 
    \alpha_W(x,\vartheta\,|\,x_0)\,\barF_J(S-x)\,{\rm d}x 
 \label{eqGX} \\
 && \hspace{1.5cm}
 +\int_0^t {\rm d}F_R(\vartheta)\int_{-\infty}^S \alpha_W(x,\vartheta\,|\,x_0)\,{\rm d}x 
 \int_{-\infty}^{S-x} G_X(S,t-\vartheta\,|\,x+u)\,{\rm d}F_J(u). 
    \nonumber
\end{eqnarray}
Hence, after setting 
\begin{eqnarray}
 && B_X(S,t\,|\,x_0):=\barF_R(t)\,G_W(S,t\,|\,x_0)
    +\int_0^t {\rm d}F_R(\vartheta)\int_{-\infty}^S 
    \alpha_W(x,\vartheta\,|\,x_0)\,\barF_J(S-x)\,{\rm d}x 
  \label{defBX} \\
 &&  \hspace{0.5cm}
 +\int_0^t \barF_R(t-\vartheta)\,{\rm d}F_R(\vartheta)\int_{-\infty}^S 
 \alpha_W(x,\vartheta\,|\,x_0)\,{\rm d}x 
 \int_{-\infty}^{S-x} G_W(S,t-\vartheta\,|\,x+u)\,{\rm d}F_J(u) 
  \nonumber 
\end{eqnarray}
it is not hard to see that from (\ref{eqGX}) the following lower bound 
for the FCT distribution function is obtained: 
\begin{equation}
 G_X(S,t\,|\,x_0)\geq B_X(S,t\,|\,x_0), \qquad t>0. 
 \label{boundG}
\end{equation}
%
\subsection{An improved bound} 
Bound (\ref{boundG}) holds in the case of generally distributed jumps $J_i$ and 
renewals $R_i$. By making use of a method already adopted in the proof of Theorem 3.3 
of Di Crescenzo and Pellerey \cite{DiCrPe}, hereafter we see how such bound can be 
improved by specifying the distributions $F_J$ and $F_R$. Let us assume that 
\par
(i) \ the random times separating consecutive jumps are exponentially distributed, 
with $\barF_R(t)=e^{-\lambda t}$, $t\geq 0$, and $\lambda>0$, and that 
\par
(ii) \ the jumps have fixed amplitude a.s., with $F_J(x)=\uno_{\{x\geq a\}}$ 
and $a>0$, which corresponds to assumption (\ref{(vaJ)}) with $b=0$ and $\eta=1$. 
\par
In the following we shall stochastically compare $\{X(t)\}$ 
with another Wiener process with jumps $\{\widetilde X_n(t)\}_{t \geq 0}$, 
which is driven by the same diffusive component of $\{X(t)\}$. We formally have
\begin{equation}
 \widetilde X_n(t)= W(t) + \widetilde Y_n(t), \qquad t \geq 0,
 \label{(mod2)}
\end{equation}
where $\{W(t)\}$ and $\{\widetilde Y_n(t)\}$ are independent stochastic processes, 
$\{W(t)\}$ is the same process appearing in the right-hand-side of (\ref{(model)}),  
$\{\widetilde Y_n(t)\}$ is a jump process such that $\widetilde Y_n(0)=0$ 
and $\widetilde Y_n(t)=\sum_{i=1}^{\widetilde N(t)} \widetilde J_i$, $t>0$, 
with $\widetilde Y_n(t)=0$ when $\widetilde N(t)=0$. Moreover,  
$\widetilde J_1,\widetilde J_2,\ldots$ are real-valued i.i.d.\ r.v.'s 
degenerating at $na$, with $a>0$ and $n$ a fixed positive integer,  
i.e.\ possessing cdf $F_{\widetilde J}(x)=\uno_{\{x\geq na\}}$. 
Furthermore, $\{\widetilde N(t)\}_{t\geq 0}$ is a counting process independent 
of $\{\widetilde J_1,\widetilde J_2,\ldots\}$ and characterized by i.i.d.\   
Erlang-distributed renewal times $\widetilde R_1,\widetilde R_2,\ldots$ 
having survival function 
$\barF_{\widetilde R}(t)=e^{-\lambda t}\sum_{j=0}^{n-1}{(\lambda t)^j \over j!}$, 
$t\geq 0$, with $\lambda>0$. We stress that processes $\{X(t)\}$ and 
$\{\widetilde X_n(t)\}$ share the two parameters $\lambda$ and $a$. 
\par
Making use of a customary tecnique based on the constructions of ``clone'' processes 
(see, for instance, Theorem 3.3 of Di~Crescenzo {\em et al.}\ \cite{DCrDNR05b}) and 
recalling (\ref{(model)}) and (\ref{(mod2)}), it is not hard to prove that 
\begin{equation}
 \widetilde X_n(t)\leq_{st} X(t) \qquad \hbox{for all } n\geq 1,
 \label{Xstocordin}
\end{equation}
where $\leq_{st}$ denotes the usual stochastic order (see Section 1.A of 
Shaked and Shanthikumar \cite{ShSh94}). An immediate consequence of (\ref{Xstocordin}) 
is that the FCT's of those processes are stochastically ordered too, i.e.
\begin{equation}
 T_X\leq_{st}T_{\widetilde X_n} \qquad \hbox{for all } n\geq 1. 
 \label{Tstocordin}
\end{equation}
In other terms, 
\begin{equation}
 G_X(S,t\,|\,x_0)\geq G_{\widetilde X_n}(S,t\,|\,x_0)
 \qquad \hbox{for all } n\geq 1  \hbox{ and }t\geq 0, 
 \label{Gordin}
\end{equation}
where $G_{\widetilde X_n}(S,t\,|\,x_0)$ is the FCT cdf of $\{\widetilde X_n(t)\}$. 
Similarly to (\ref{boundG}) the following bound holds for all $n\geq 1$: 
\begin{equation}
 G_{\widetilde X_n}(S,t\,|\,x_0)\geq B_{\widetilde X_n}(S,t\,|\,x_0), \qquad t>0, 
 \label{boundtG}
\end{equation}
where 
\begin{eqnarray}
 B_{\widetilde X_n}(S,t\,|\,x_0) 
 \!\!\! &:=& \!\!\! e^{-\lambda t}\sum_{j=0}^{n-1}{(\lambda t)^j\over j!}\,G_W(S,t\,|\,x_0)
 \label{integrale} \\
 &+& \!\!\! \int_0^t \lambda e^{-\lambda \vartheta}{(\lambda\vartheta)^{n-1}\over (n-1)!}
 {\rm d}\vartheta\int_{S-na}^S \alpha_W(x,\vartheta\,|\,x_0)\,{\rm d}x 
 \nonumber \\
 &+& \!\!\!
 \int_0^t e^{-\lambda (t-\vartheta)}\sum_{j=0}^{n-1}{[\lambda(t-\vartheta)]^j \over j!}
 \lambda e^{-\lambda \vartheta}{(\lambda \vartheta)^{n-1} \over (n-1)!}{\rm d}\vartheta 
 \nonumber \\
 &\times & \!\!\!
 \int_{-\infty}^{S-na} \alpha_W(x,\vartheta\,|\,x_0)\,
 G_W(S,t-\vartheta\,|\,x+na)\,{\rm d}x, 
 \nonumber 
\end{eqnarray}
with $G_W$ and $\alpha_W$ defined in (\ref{(distr)}) and (\ref{(avoid)}), respectively.
In conclusion, by (\ref{Gordin}) and (\ref{boundtG}) we obtain 
\begin{equation}
 G_{X}(S,t\,|\,x_0)\geq 
 \sup_{n\geq 1} \,B_{\widetilde X_n}(S,t\,|\,x_0), 
 \qquad t\geq 0. 
 \label{disfinale}
\end{equation}
Since (\ref{boundG}) provides a bound of type 
$G_{X}(S,t\,|\,x_0)\geq B_{\widetilde X_1}(S,t\,|\,x_0)$,   
Eq.\ (\ref{disfinale}) yields a better bound. Furthermore, the 
right-hand-side of (\ref{disfinale}) is not necessarily a distribution 
function, so that we can improve the bound as follows: 
\begin{equation}
 G_{X}(S,t\,|\,x_0)\geq G_{\ell}(S,t\,|\,x_0)\equiv
 \max_{0\leq \tau\leq t}\,\sup_{n\geq 1} \,B_{\widetilde X_n}(S,\tau\,|\,x_0), 
 \qquad t\geq 0. 
 \label{disultima}
\end{equation}
%
\begin{figure}[t]
\begin{center}
\leavevmode
\epsfxsize=300pt
\epsffile{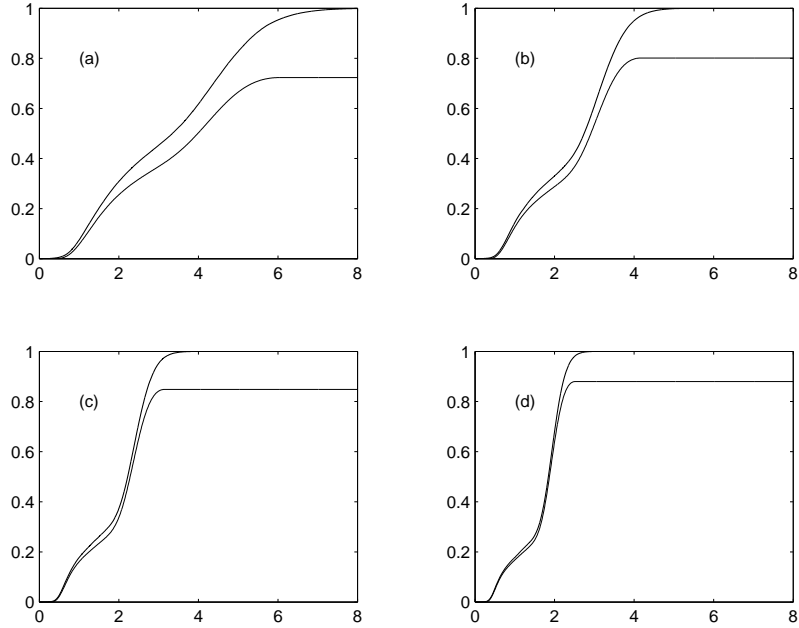}
\caption{
Estimated FCT distribution function $\widehat{G}_X(S,t\,|\,x_0)$ and lower bound 
$G_{\ell}(S,t\,|\,x_0)$ for the same cases treated in Figure 1. 
}
\end{center}
\end{figure}
%
\begin{figure}[t]
\begin{center}
\leavevmode
\epsfxsize=180pt
\epsffile{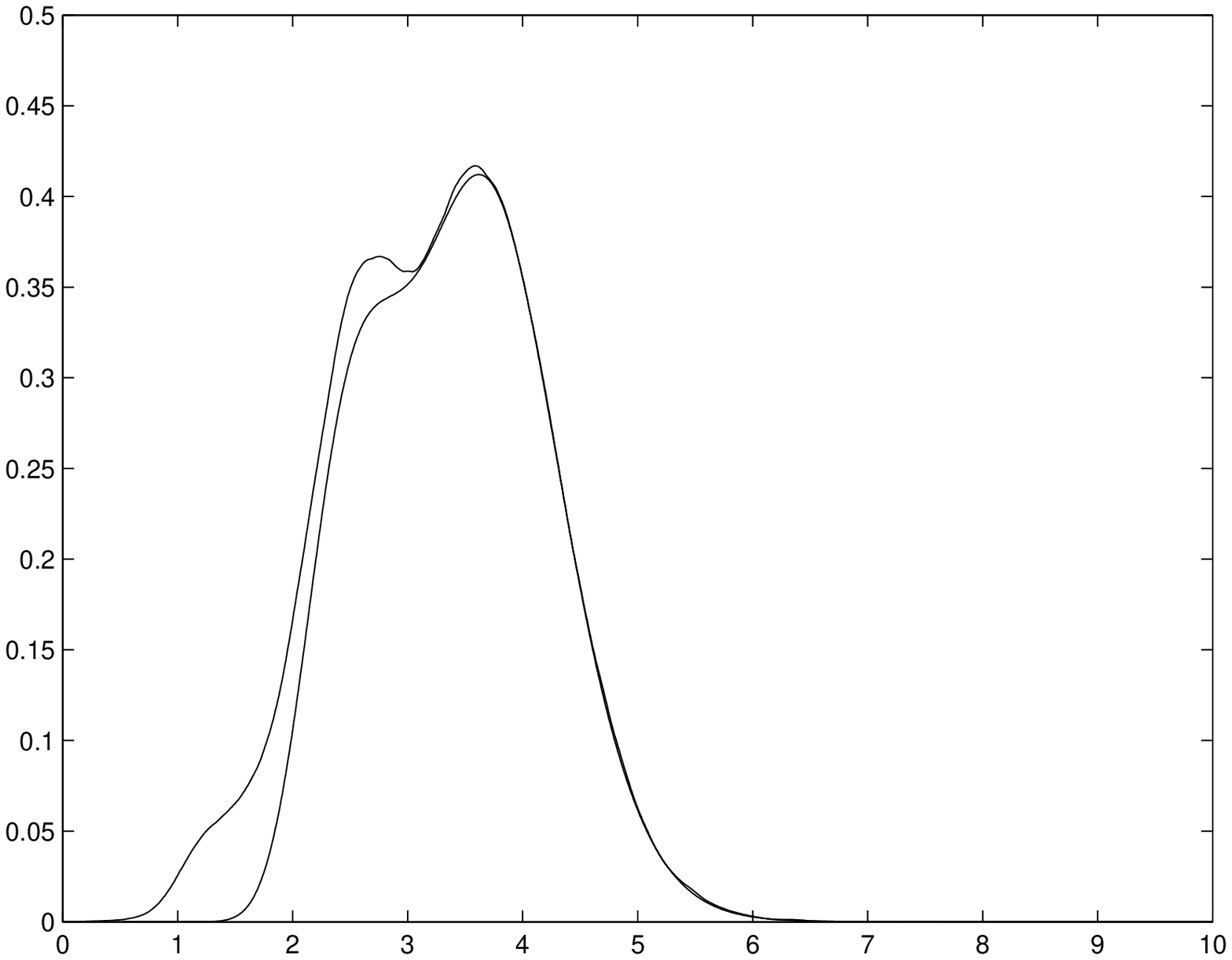}
\qquad
\epsfxsize=180pt
\epsffile{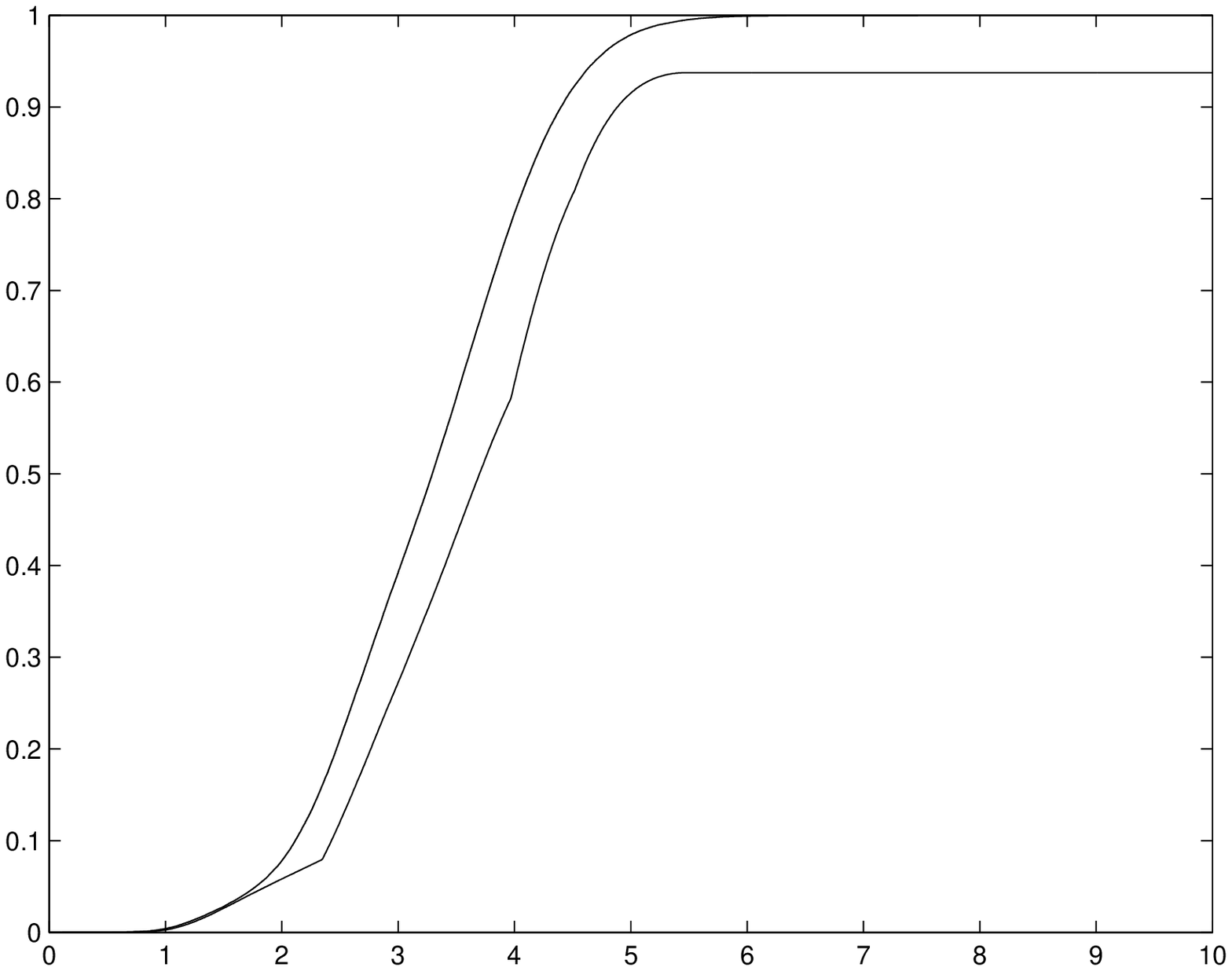}
\caption{
On the left: Estimated FCT density $\widehat{g}_X(S,t\,|\,x_0)$ and lower bound 
$g_{\ell}(S,t\,|\,x_0)$ for $a=2$, $b=0$, $x_0=0$, $S=6$, $\lambda=0.2$, $\eta=1$, $\mu=1.5$ 
and $\sigma^2=0.2$. On the right: The corresponding estimated distribution function 
$\widehat{G}_X(S,t\,|\,x_0)$ and lower bound $G_{\ell}(S,t\,|\,x_0)$. 
}
\end{center}
\end{figure} 
\par
Figure 3 shows estimates of FCT distribution function 
$\widehat{G}_X(S,t\,|\,x_0)$ obtained from simulations perfomed by means 
of the procedure given in Di~Crescenzo {\em et al.}\ \cite{DCDNR2005}, 
together with the respective lower bound (\ref{disultima}). We stress that 
the lower bounds given in Figure 3 shows a case in which, for all $\tau\geq 0$ 
\begin{equation}
 \sup_{n\geq 1} \,B_{\widetilde X_n}(S,\tau\,|\,x_0)=B_{\widetilde X_1}(S,\tau\,|\,x_0).
 \label{fsup}
\end{equation} 
Note that (\ref{fsup}) holds because relation $S-na>0$ is satisfied only if $n=1$. 
Instead, the angular points appearing in the lower bound on the right of Figure 4   
show a case in which 
$$
 \sup_{n\geq 1} \,B_{\widetilde X_n}(S,\tau\,|\,x_0)=B_{\widetilde X_k}(S,\tau\,|\,x_0),
$$ 
for different values of $k$ as $\tau$ varies. 
%
\section{Remarks on computational aspects} 
Let us now point out some computational features related to the evaluation 
of bound (\ref{disultima}). The last integral in (\ref{integrale}), i.e. 
$$
 \int_{-\infty}^{S-na} \alpha_W(x,\vartheta\,|\,x_0)\,
 G_W(S,t-\vartheta\,|\,x+na)\,{\rm d}x, 
$$ 
can be expressed in terms of integrals of the form
\begin{equation}
 I = \int_{0}^{+\infty} e^{-a(x-b)^2} \,\Phi(cx+d)\,{\rm d}x, 
 \label{(intnum)}
\end{equation}
with $a > 0$, with $b,c,d \in \Reali$ and $\Phi$ given in (\ref{(norm)}). 
Unfortunately no closed forms are available for (\ref{(intnum)}), so that 
we have been forced to numerically evaluate it by splitting the
integration domain as follows:
$$
 I = \int_{0}^{b} e^{-a(x-b)^2} \,\Phi (cx+d)\,{\rm d}x 
 + \int_{b}^{+\infty} e^{-a(x-b)^2} \,\Phi (cx+d)\,{\rm d}x.
$$ 
Then, transforming the second integral we obtain:
\begin{equation}
 I=\int_{0}^{b} e^{-a(x-b)^2} \,\Phi (cx+d)\,{\rm d}x
 +\frac{1}{2 \sqrt{a}}\int_0^{+\infty}x^{1/2} e^{-x}\,
 \Phi\left(c\sqrt{\frac{x}{a}} + c \, b + d \right) {\rm d}x.
 \label{ultint}
\end{equation}
A Gauss-Legendre quadrature rule has been applied to the first integral by choosing 
$16$ abscissas, after an {\em ad hoc\/} transformation of the integration domain 
in $[-1,1]$, whereas a generalized Gauss-Laguerre quadrature rule has been applied 
to the second integral by choosing $16$ abscissas. Note that the latter quadrature 
formula is exact for weighting functions as $w(x)=x^{\beta} e^{-x}$, with 
$\beta\in\Reali$, which includes the case of the second integral in (\ref{ultint}). 
\begin{table}[t]
\begin{center}
\begin{tabular}{c|ccccc}
\hline
$\mu\;\backslash\;\eta$ & 1	& 0.8	& 0.5	& 0.2	& 0.0 \\
\hline
1.0	& 3.9358e-02	& 1.1393e-01	& 2.4054e-01	& 3.6791e-01	& 4.5211e-01 \\
1.5	& 2.1059e-02	& 6.5103e-02	& 1.4572e-01	& 2.2688e-01	& 2.8409e-01 \\
2.0	& 1.7266e-02	& 4.4515e-02	& 9.6393e-02	& 1.5005e-01	& 1.9656e-01 \\
2.5	& 1.7283e-02	& 3.4975e-02	& 7.0264e-02	& 1.1009e-01	& 1.4486e-01 \\
\hline
\end{tabular}
\caption{Values of ${\cal E}$ for various choices of $\mu$ and $\eta$.}
\end{center}
\end{table}
\par
In order to evaluate the goodness of lower bound (\ref{(bound)})
we shall now adopt the following ``measure of closeness'': 
\begin{equation}
 {\cal E}=h\sum_{j=1}^n \left[\widehat{g}_X(S,jh\,|\,x_0)-g_{\ell}(S,jh\,|\,x_0)\right],
 \label{eq:errore}
\end{equation}
where $h$ is a discretization step and $n$ is the minimum integer such that 
$$
 h\sum_{j=1}^n \widehat{g}_X(S,jh\,|\,x_0) \geq 0.99.
$$
Differently stated, (\ref{eq:errore}) measures how the bound $g_{\ell}$ is close to 
$\widehat{g}_X$ by the multiples of $h$ up to the $0.99$-percentile. For $h=0.01$, 
Table 1 shows values of ${\cal E}$, where $\widehat{g}_X$ is an estimate of $g_X$, 
whereas $g_{\ell}$ is the lower bound given in (\ref{(bound)}). It is evident that 
${\cal E}$ increases if $\mu$ decreases or if $\eta$ decreases. 
\par
We finally point out that throughout the paper the estimated functions 
$\widehat{g}_X$ and $\widehat{G}_X$ have been obtained by use of $10^6$ simulated FCT's. 
The estimated pdf $\widehat{g}_X$ has been built by adopting an Epanechnikov 
kernel estimator, while $\widehat{G}_X$ derives by use of a Kaplan-Meier estimator. 
%
\subsection*{\bf Acknowledgments}
This work has been performed under partial support by MIUR (PRIN 2005), 
by G.N.C.S.-INdAM and by Campania Region.

%
%
{\em Antonio Di Crescenzo}: Dipartimento di Matematica e Informatica, 
Universit\`a degli Studi di Salerno, Via Ponte don Melillo, 
I-84084 Fisciano (SA), Italy. E-mail: adicrescenzo@unisa.it \\ 
{\em Elvira Di Nardo}: Dipartimento di Matematica e Informatica, 
Universit\`a degli Studi della Basilicata, 
Campus Macchia Romana, I-85100 Potenza, Italy. E-mail: dinardo@unibas.it \\
{\em Luigi M.\ Ricciardi}: Dipartimento di Matematica e Applicazioni, 
Universit\`a di Napoli Fe\-derico II, 
Via Cintia, I-80126 Napoli, Italy. E-mail: luigi.ricciardi@unina.it
%
\end{document}